\documentclass{amsart}
\usepackage{amssymb}
\newtheorem{theo}{Theorem}[section]
\newtheorem{lemm}[theo]{Lemma}

\newtheorem{step}{Step}

\newcommand{\E}{\mathcal{E}}
\newcommand{\cb}{\mathcal{B}}

\newcommand{\cg}{\mathcal{G}}

\begin{document}

\title{Planar Graphs and Covers}
\author{M.J. Dunwoody}

\maketitle

Version 28.4.2009

\begin{abstract}

Planar locally finite graphs which are almost vertex transitive are discussed. If the graph is $3$-connected and has at most one end then the group of automorphisms is a planar discontinuous group and its structure is well-known. A general result is obtained for such graphs where no restriction is put on the number of ends. It is shown that such a graph can be built up from one-ended or finite planar graphs in a precise way. The results give a classification of the finitely generated groups with planar Cayley graphs.

\end{abstract}

 \section{ Introduction}
 
 In 1965 Maskit \cite {[Ma]} proved his Planarity Theorem which classified the planar covers of compact surfaces.  
 This classification is significant in understanding planar graphs. Thus one way
 of obtaining a planar graph is by embedding a finite graph $Y$  in a surface and then considering
 its lift in a regular planar cover.  A planar graph $X$ obtained in this way will have an  automorphism
 group $G$ such that $G\backslash X$ is finite. 
 We show any $3$-connected planar graph $X$ with an automorphism group $G$ such that 
 $Y = G\backslash X$ is finite arises in this way. Maskit's Theorem enables one to classify
 the groups $G$ that arise and we show that each such group has a planar Cayley graph.
 Thus one has a classification of the groups that admit a $3$-connected, locally finite, planar Cayley graph.
 
 In Section~\ref{section:2conn} graphs that are not $3$-connected are discussed.
 It is shown that if $X$ is such a graph which  is connected and  locally finite with automorphism group $G$ then
 there is a $G$-tree $T$, associated with $X$.
 To each vertex $v \in VT$ with stabilizer $G_v$, either $G_v$ contains a normal subgroup $N_v$ which fixes a vertex of $X$ and $G_v/N_v$ is either a finite cyclic or dihedral group   or there is  a $3$-connected $G_v$-graph $Z_v$ which
 is a subdivision of a subgraph of $X$ and
 is planar if $X$ is planar.  Each edge $e$ of $T$ is associated with a decomposition 
 $X = A\cup A^*$ where $A, A^*$ are subgraphs $A\cup A^* = X$ and $A\cap A^*$ has either one or two vertices and no edges.   This is a new proof of a result of Droms, Servatius and Servatius \cite{[DrSS]} which generalizes the analysis of Tutte \cite{[Tu]} of finite $2$-connected graphs.
 A consequence of this result is that the $3$-connected condition can be removed from
 our classification of locally finite, planar Cayley graphs.
 The theory here is similar to the construction of  trees in \cite{[DDu]}, Chapter II.
 The difference is that in \cite{[DDu]} we use finite disconnecting edge sets rather than disconnecting vertex sets.  Thus in \cite{[DDu]} we consider subsets of $VX$ which have finite coboundary.  \footnote{ Recently (May, 2009) in \cite{[DK]},   Bernhard Kr\" on and I obtain an even more general result further analyzing the situation when
 a graph can be disconected by removing a finite set of vertices.}
 
 The coboundary of the set $A \subseteq VX$ is the set 
 $$\delta A = \{  e\in EX | \ e \ {\rm joins\  vertices
 \ in} \ A \  {\rm and}  \ VX - A \}.$$
   The set of all such subsets is denoted ${\cb}X$, and is a Boolean
 ring, as it is closed under addition, multiplication, and complementation.
 The graph $X$ is defined to be {\it accessible} if there is an integer $k$ such that ${\cb}X$ is generated (as a ring) by the set $A \in {\cb}X$ such that $|\delta A| \leq k$.
 In \cite{[TW]} it is shown that this is equivalent to the condition that any two ends of $X$ can be separated
 by removing at most $k$ edges.  In \cite{[DDu]} it is shown that if $X$ is accessible then
 ${\cb}X$ is generated by a tree set, i.e. a set in which the partial order induced by subset inclusion
 can also be interpreted as the natural order on the directed edges of a tree.
 If a group $G$  is accessible, i.e. if it has an accessible Cayley graph, then it has 
 a decomposition as the fundamental group of graph of groups in which edge groups are finite and
 vertex groups are finite or have one end.

 In \cite{[Du]} I showed that the Cayley graphs of finitely presented groups are accessible.  This result could be viewed
 as a generalization of Maskit's result, though at the time I was unaware of that fact.
 Later in \cite{[DDu]} an account is given of my earlier result and also how it applies in the situation
 of regular planar covers.  Thus a proof is given there of Maskit's Theorem without 
 acknowledgment of his work.  I apologize for this omission.
 In \cite{[Du2]} I gave an example of a finitely generated inaccessible group.
 Thus not every locally finite vertex transitive graph is accessible.  In later papers I gave further
 examples of such graphs.  
 
 In this paper it is shown that if $X$ is a locally finite planar graph with automorphism group
 $G$, and $G\backslash X$ is finite, then $X$ is accessible.  This then enables one to
 give very precise information about the structure of both $X$ and $G$.  In particular
 it is possible to classify the groups which have planar Cayley graphs.
 
 Planar Cayley graphs have been studied in the papers of C. Droms and C. Droms, B. Servatius
 and H. Servatius.  Thus in \cite{[Dr]} it is proved that a plane group (i.e. a finitely generated group
 with a plane Cayley graph) is accessible and in \cite{[DrSS2]}  a result is proved that is very close to
 the classification obtained here, except that Maskit's  Planarity Theorem is not applied to classify
 planar covers.

  I thank Gareth Jones for showing me the paper \cite{[Mo]} and for helpful discussions.

 \bigskip
 \noindent
 \section{Planar graphs}

  A facial path in a plane graph is one obtained by taking a directed  edge and then proceeding in the direction of the edge
 always turning sharp left and in the reverse  direction always turning sharp right.  One then obtains either a cycle or a $(2, \infty )$-path.
Clearly every directed edge uniquely determines such a path, and so each unoriented edge lies in
 two such paths.  
 Recall that a graph is $3${\it -connected} if it has at least $5$ vertices and no set of $2$ vertices disconnects it.  By Menger's Theorem this is equivalent to saying that a graph is $3$-connected if it has 
 at least $2$ vertices and every pair of distinct vertices is joined by at least $3$ internally disjoint paths.
 Let $X$ be a $3$-connected planar graph. It was proved by Whitney \cite{[W]} for finite graphs that an embedding of $X$ in the $2$-sphere  is essentially unique, and Imrich showed that for an infinite graph the cyclic order of edges at a vertex is uniquely determined.  In \cite{[RT]} it is proved that the Freudenthal compactification of a planar graph has a unique embedding in the $2$-sphere.
Thus we have

\begin{theo}\label{theo:facial}  Let $X$ be a $3$-connected locally finite planar graph.  Then any automorphism of $X$ takes facial paths to facial paths.
     \end{theo}

  Let $X$ be  a locally finite planar graph $3$-connected graph with automorphism group $G$ 
 so that $G\backslash X$ is a finite graph $Y$.  
It follows from Theorem~\ref{theo:facial} that we
 can attach $2$-cells to the finite cycles in $X$ which bound faces in the plane graph
 and obtain a cellular $2$-complex $K$ with polyhedron $M = |K|$ which admits an action of $G$ so that $G\backslash M$ is compact.  It may be the case that a vertex can lie in more than one infinite facial path (see \cite{[Mo]}),
 in which case $M$ is not a $2$-manifold with boundary.  However we can embed
 $X$ in a graph $\hat X$ which admits an action of $G$ and for which the corresponding
 complex $\hat M$ is a $2$-manifold with boundary.  This is done by identifying  each
 infinite facial path with  one of the two sides  of an infinite ladder. That one can do this follows from the previously
 quoted result of \cite{[RT]}.
 Thus we can assume that $M$ is a planar $2$-manifold with boundary and that 
 $G\backslash M$ is compact.
 We also assume that the action of $G$ on $M$ is orientation preserving as if this is not the
 case then $G$ has a subgroup of $H$ index $2$ for which the action is orientation preserving
 and $H\backslash M$ is compact.

 Our main tool in the classification is the theory of tracks in $2$-complexes (see  \cite{[DDu]}, Chapter
 VI).  The argument of \cite{[DDu]} is generalized to that of compact orbifolds  with boundary as this
  provides the most effective application to planar graphs.

 We reproduce p.245 of {\cite{[DDu]} with some corrections and alterations.
 The numbering of Theorems and Propositions is that of \cite{[DDu]}.
 As noted in \cite{[DDu]}, Errata the condition that $H^1(K. {\bf Z}_2) = 0$ is much too strong, and rules out
 most of the planar spaces that it was intended to study.  Instead one just needs that
 simple closed curves (scc's) separate.
 
 Let, then,  $M = |K|$ be a $2$-manifold possibly with boundary  and in which
 every track corresponding to a scc separates.   Here we are assuming that $K$ is a simplicial
 complex though later we will only assume that it is a cell complex. Suppose there is a group $G$ which acts on $K$ and which thereby induces an
 action on $M$. If we assume that $G$ acts faithfully on $K$  (so that it is a subgroup of automorphism
 group of the $1$-skeleton) then the stabilizer of each simplex $\sigma \in SK$ must be finite.
 This is because any automorphism which fixes a vertex $v_0$ and its incident edges must fix the whole 
 graph.  This follows by considering a sequence of finite subcomplexes $\{ v_0 \} = K_0 \subset K_1
 \subset K_2 \subset \dots $ whose union  is $K$ and such that $K_i$ is obtained from $K_{i-1}$ by
 adjoining an edge at least one vertex of which is in $K_{i-1}$, and showing that the automorphism
 must fix each term of the sequence.
 If we further assume that the action is orientation preserving -  either $G$ or a subgroup of index $2$
 will have this property -
 then for each $v \in VK, G_v$ must be a finite cyclic group and the stabilizer of an edge is trivial.

 We also know that $G\backslash K = L$ is finite, since $G\backslash M$ is compact.
 Suppose there is an scc in $M$ which does not bound a disc.  By Proposition 7.1 there is
 a track which is an scc that does not bound a disc.
 Such a track is called {\it thin} if it has the minimal number of intersections with the one
 skeleton of $K$. Note that for the moment we do not consider tracks which
 are arcs connecting boundary points of $M$.
 Let $t_1, t_2$ be thin tracks. In this special situation we want to prove an analogous result
 to that of Proposition 5.5, i.e that $t_1 + t_2 = s_1 \bigoplus s_2$ there $s_1$ and $s_2$
 are disjoint thin tracks. Note that a thin track separates $M$ into two components each
 containing infinitely many vertices.  For if one of the components contained only finitely
 many vertices then it has empty intersection with $\partial M$ since each component of
 $\partial M$ is infinite.  The closure of the component containing finitely many vertices will be a compact planar $2$-manifold with boundary
 $S^1$ which must be a disc. Conversely,  if a track separates $M$ into two components
 with infinitely many vertices neither can be a disc.
 The argument of Proposition 5.5 now goes through.  As is noted in the \cite{[DDu]}, Errata
 after line 9 of p236 the following sentence should be added -  ``Moreover it follows from the thinness of
 $b_2^*$ or $b_1^*$ that $\nu = \delta$." 
 By Theorem 5.9 there is an scc $\ell $ such that $g\ell =\ell$ or $g\ell \cap \ell =\emptyset $
 for every $g \in G$ and $\ell $ does not bound a closed disc in $M$. In fact we could
 take $\ell $ to be a thin track as in Section 5.
 
 Let $P$ be a $G$-pattern in $M$, where the component tracks of $P$ are scc's.  Let $C$ be the closure in $M$ of a component of $M-P$. Thus $C$ is a two-manifold (with boundary) acted on by $G_C$, and $\partial C$
 contains  a $G_C$-set of scc's. We obtain a two-manifold $\hat C$ (possibly with boundary components which are copies of $\bf R $)  by attaching discs to each
 boundary scc in $\partial C$.  If $\hat C$ has an scc $\ell $ which does not bound a disc in $\hat C$ 
 then there is such a curve which actually lies in $C$.  Again as in the proof of Theorem 5.9,
 we can choose $\ell $ so that it is a track, disjoint from $P$ and either $g\ell = \ell$ or 
 $g\ell \cap \ell = \emptyset $ for every $g \in G$.  Thus we can replace $P$ by $P\cup
 \{ g\ell | g\in G \} $.  Since $G\backslash K = L$ is finite, we cannot repeat this process more
 than $n(L)$ times , and so we eventually find a $G$-pattern $P$ in $M$, such that, if 
 $\hat C $ is obtained by attaching discs to each boundary scc in the  closure of a component of $M-P$ then every scc in $\hat C$ bounds a disc. If $\hat C = S^2$ then its stabilizer is a finite group of 
 isometries of $S^2$ while
  if $\hat C = \bf R^2$ then its stabilizer admits an action by isometries
 on either the Euclidean or hyperbolic plane.  In both cases $\hat C$ has empty boundary
 and it will be a planar discontinuous group.  The structure of
 these groups is well known (see \cite{ [ZVC]}). 
 
 If $\hat C$ is not $\bf R^2$ or $S^2$  then it must contain at least one component of $\partial M$.
 In fact $\hat C $ must contain more that one component.  For if not then $G_C$ is the stabilizer of
 the single boundary component $t$ which is infinite cyclic and there is no cocompact action of an infinite cyclic group on a simply connected two-manifold with a single boundary component.
 Consider an arc $\ell $  joining two boundary components.  Then $\ell $ must separate $C$.
 For if not there is an scc in $C$ meeting $\ell $ in a single point.  But this curve
 will bound a disc in $\hat C$, and  this disc must contain a component of $\partial M$, since it contains an end point of $\ell$.
But each component of $\partial M$ is not compact and so we have a contradiction. Hence $\ell $ separates.
 Since every arc joining points in $\partial \hat C$ separates, we can repeat the 
 process described in the previous paragraph for scc's, to obtain a $G_C$-pattern $P_C$ of tracks which are arcs   such that the closure of each component of 
$\hat C - P_C$ is a disc.  
Each such disc will have finite cyclic stabilizer.  Note that the stabilizer of each arc in $P_C$ is trivial or cyclic of order two
and each track with non-trivial stabilizer can be replaced by a parallel pair for which
the stabilizers are trivial.  Thus  $G_C$ acts
on a tree with trivial edge groups and finite cyclic vertex groups.
By Bass-Serre theory a group has such an action if and only if it is a free product of cyclic groups
(not necessarily finite). Such a group will be the fundamental group of a graph of finite cyclic groups.  We have proved the following.

\begin{theo}\label{theo:main} If  $G$ is the automorphism group of a $3$-connected locally finite planar graph $X$ with
$G\backslash X$ finite, then $G$ (or an index two subgroup of $G$) is the fundamental group
of a graph of groups in which each vertex group is either a planar discontinuous group
or a free product of finitely many  cyclic groups, and all edge groups are finite cyclic groups (possibly
trivial).  
\end{theo}

This means that any such group can be built up from planar discontinuous groups and finite cyclic groups
in finitely many steps by free products with amalgamation along a finite cyclic subgroup
or by HNN extensions along a finite cyclic subgroup. This sort of result was conjectured
in \cite{[Mo]}.

In fact our proof gives more information than this as it provides information about the action
of $G$ on the $2$-manifold $M$.
Let ${\cg}(Y)$ be the graph of groups in the above, so that for each
vertex $v \in VY$ the group $G(v)$ is either planar discontinuous  or a free product of
finitely many cyclic groups, and for each oriented edge $e \in EY$ the edge group $G(e)$
is finite cyclic and there are specified injective homomorphism $G(e) \rightarrow G(\iota e),
G(e) \rightarrow G(\tau e)$ into the associated vertex groups.

The additional information is contained in the description of the orbifolds  $M$ and $G\backslash M = W$.
Thus for each vertex $v \in VY$ there is a compact $2$-orbifold  $W(v)$ which is a compact
surface - possibly with boundary - and finitely many special interior points each with a specified index 
which is a positive integer. For each edge $e \in EY$ there is an arc joining $a(e)$ joining
a special point in $W(\iota e)$ and a special point in  $W(\tau e)$ with the same index.
We thicken each arc so that it it is the axis of a solid cylinder whose end discs are small closed discs
about the special points in $W(\iota e)$ and $W(\tau e)$.  Now we can take $W$ to be the
union of the $W(v)$'s minus the interior of any disc that is the end of an edge cylinder together
with the closed annulus that is the curved part of the boundary of each such cylinder.
Note that each special point cannot be the end point of more than one $a(e)$ but there may
be special points which are not the end point of any arc.  For these points it can be assumed the
index is more than one.
There is a well-known classification of closed surfaces and closed surfaces with boundary
and also of compact $2$-orbifolds (see \cite{[Sc]}), and so can classify the spaces $W$.
Note that $W$ is itself a $2$-orbifold in which the special points are those in the vertex 
$W(v)$'s that are not the end points of any $a(e)$. There is  also a set $P$ of scc's in $W$,
each element of $W$
corresponding to  a meridian in the annulus  which is  the curved part of the boundary  of a cylinder around a particular $a(e)$.
Each such scc has a $k$-fold covering by an scc in  $M$ where $k$ is the index shared by the
two end points of the corresponding arc $a(e)$.

It can be seen from the above that any group $G$ with the structure given in the 
theorem can be realized as the automorphism group of a planar graph.
In fact for an appropriately chosen set of generators $S$ the Cayley  graph of $G$ with respect
to $S$ will be planar.
The group $G$ has a presentation \medskip
$$ G = \langle a_1, b_1, a_2, b_2, \dots a_p, b_p, e_1, e_2, \dots , e_r , f_1, f_2 , \dots f_s |
e_1^{m_1}=  e_2^{m_2} \dots  = e_r^{m_r} $$
$$= w_1^{n_1} = w_2^{n_2} \dots = w_t^{n_t} = [a_1, b_1][a_2, b_2]\dots [a_p, b_p]e_1,e_2\dots e_rf_1f_2\dots f_s = 1 \rangle , $$
where $m_1, m_2, \dots , m_r$ are integers bigger than one, $w_1, w_2, \dots , w_t$ are words
in the generators representing the scc's in $P$ and the positive integers $n_1, n_2, \dots , n_t$
are the corresponding values of the index $k$.

Consider the group $\bar G$ with presentation
$$ \bar G = \langle a_1, b_1, a_2, b_2, \dots a_p, b_p, e_1, e_2, \dots , e_r , f_1, f_2 , \dots f_s | $$
$$[a_1, b_1][a_2, b_2]\dots[a_p, b_p]e_1e_2\dots e_rf_1f_2\dots f_s = 1 \rangle , $$
The Cayley complex $\bar W$ corresponding to this presentation  has $2p + r + s $ loops corresponding to the
generators all sharing a single base point $b$ to make a $(2p+r + s)$-leaf rose $R$ and there is a single $2$-cell attached along the
relator word.  The universal cover $\bar M$  of $\bar W$ is planar and the lift of $R$ is the Cayley graph $\bar X$  of $\bar G$ with respect to the set of generators.
Let $N$ be the normal subgroup of $\bar G$ which is the normal closure of the elements corresponding
to the other relators of $G$.  Then $M = N\backslash \bar M$ is planar and there is an
embedding of the Cayley graph $X = N \backslash \bar X$ of $G$ in $M$.

What we have proved here is essentially a generalization of Maskit's Planarity Theorem
to surfaces with boundary.

 \section{Graphs which are not $3$-connected}\label{section:2conn}
 \medskip
 \def\lq{\subseteq}
  Let $X$ be a locally finite connected graph which is $1$-connected but not $2$-connected.
 This means that there are cut-points, i.e. vertices whose removal disconnects the graph.
 
 We consider subgraphs  $A \subset X$ with the following property
 \begin{itemize}
 \item [{(i)}]  there is exactly one vertex $u \in VA$ such that if $w \in VA, w \not= u$
 then every edge of $X$ incident with $w$ is also in $A$.  Some, but not all, edges incident with
 $u$ are in $A$.
 \end{itemize}
 Define $\delta A = \{ u \}$.
 
 Let ${\cb}_1X$ denote the set of all subgraphs of $X$ satisfying (i).
 Let $A \in {\cb}_1X$.  Let $A^*$ be the subgraph
 with $VA^* = (VX - VA)\cup \{ u\} $ and $EA^* = EX - EA.$
Clearly $A^*$ satisfies (i) and so $A^* \in {\cb}_1X$ with $\delta A^* = \{ u \}$.
Also if $G = Aut X$ so that $X$ is a $G$-graph (acting on the left) then
${\cb}_1X$ is invariant under $G$.

Since $X$ is locally finite removing a cut-point $x$ produces finitely many
connected components $C_1, C_2, \dots , C_k$  the union of any proper subset of
the $C_i$'s together with the edges joining them to $x$ and also including $x$, gives an
element of ${\cb}_1X$ and any element of ${\cb}_1X$ arises in this way.  Let $\E$ be the subset of ${\cb}_1X$ containing only those elements $A$ for which just one of $A$ or $A^*$ corresponds to a single component
$C_i$.  If $A$ corresponds to one component then it is a  minimal element of ${\cb}_1X$ containing 
$x$.  If both $A$ and $A^*$ correspond to a single component $C_i$, then $k=2$ and $x$ is a cut-point
of degree $2$,  and we can replace $X$
by a homeomorphic graph in which $x$ no longer occurs. 
A graph in which there are cut-points and in which every cut point has degree $2$ must be homeomorphic to an interval.  Thus, unless this is the case, $\E$ will be non-empty if $X$ is not $2$-connected.

\begin{lemm}\label{lemm:nesting}  Let $A, B \in \E$,
then at least one of
$A \subseteq B,  A \subseteq  B^*, A^* \subseteq B, A^* \subseteq B^*$.
\end{lemm}
\noindent
\begin{proof}  
Let $\delta A = \{ x \}$, and let $\delta B = \{ y \}$ and suppose 
that $A, B$ are minimal elements of ${\cb}_1X$ containing $x, y$ respectively, i.e. they correspond
to single components of $X - \{ x \}  , X - \{ y \} $ respectively.
If $x = y$ then either $A = B$ or $A\cap B = \{x \}$ so that $ A \lq B^*$.
Suppose $x \not= y$.
If $y \in A$ and $x \in B$ consider $A^*\cap B^*$, it is clear that any vertex in this subgraph
must have all its incident edges also in the subgraph, which, since $X$ is connected,
must be either $X$ or $\emptyset $.  It cannot be $X$ and so it is $\emptyset $ and so
$A^* \lq B$. All the other cases are treated similarly.

  Now $\lq $ is a partial order on subgraphs and $A \lq B$ implies $B^* \lq A^*$.
  Also if $A, B \in \E$ and $A\lq B$ then there are only finitely many $C \in \E$ such
  that $A \lq C \lq B$.  To see this, let $u \in VA, w \in VB^*$ and 
  let $u =v_1, v_2, \dots , v_n = w$ be the vertices of a path joining them.
  If $A \lq C \lq B$ then $\delta C = \{v_i \}$ for some $i = 1, 2, \dots , n$.
  But there are only finitely many elements of $\E$ corresponding to  a given cut-point. 
  \end{proof}

We see then that   the set $\E, \lq$ is a partially ordered set closed under the order reversing
involution $A \mapsto A^*$ and in addition it satisfies  the two conditions
\begin{itemize}
\item [{(i)}] for every $A, B \in \E$ at least  one of $A \lq B, A\lq B^*, B \lq A, B \lq A^*$ holds.
\item [{(ii)}] if $A, B \in \E$ and $A \lq B$ then there are only finitely many $C \in E$ such
that $A \lq C \lq B$.  
\end{itemize}
Note also the partial order is preserved under the action of $G$.
A $G$-set $(\E , \lq )$ satisfying the above conditions can be regarded as the oriented edge set of a $G$-tree $T = T(\E )$ (see\cite{[Du]} or (\cite{[DDu]}, p50]). Let $v \in VT$.
Different ways of defining the a vertex of $T$ are given in \cite{[Du]} and \cite{[DDu]}.  In 
\cite{[DDu]} it is defined as a particular orientation of the edge set, intuitively it is the orientation in which
the edges point towards the specified vertex.  in \cite{[Du]} a vertex is defined as the set of oriented edges
which have the specified vertex as initial vertex.
Thus one defines an equivalence relation $\sim $ on $\E $ as follows:-

$A\sim B$ if  $A \subset  B^*$ but for no $C \in \E$ is $A\subset C \subset B^*$.

It is easy to show (see \cite{ [Du]}) that  $(\E , \sim )$ is an equivalence relation and 
the natural map $\iota : \E \rightarrow VT = \E / \sim $ determines the tree $T$.
We follow this approach, so that $v \in VT$ is an equivalence class of elements of 
$\E $.

The structure of $T$ is set out in the following.

\begin{theo}\label{theo:2coonn}  Let $X$ be a connected $G$-graph, with automorphism group $G$.
There is a $G$-tree $T$ with the following properties.  The vertex set $VT$ is partitioned
$VT = J \cup K$ so that each edge $e \in ET$ has one vertex in $J$ and one in $K$.
The vertices of $J$ can be regarded as the set of cut-points of $X$ and each $u \in K$ corresponds 
to a maximal $2$-connected subgraph $B$ (or $2$-block).
Every 
vertex of $X$ lies in at least one $2$-block, and it lies in more than one if and only if it is 
a cut-point. Every edge of $X$ lies in exactly one $2$-block.  A vertex in $K$ may have infinite degree, but the vertices of $J$ have finite degree.
The $2$-block $B$ associated to the vertex $u \in K$ is a $G_u$-graph.  
\end{theo}
\begin{proof}  Most of the properties listed follow from the construction of $T$.  
If $c$ is a cut-point of $X$, then the elements of $\E $ associated with $c$,
i.e. corresponding to a single component of $X - \{ c \}$ satisfy the conditions for a vertex of
$T$. However there are other vertices of $T$.
Let $B$ be a maximal $2$-connected subgraph.  For every $A \in \E$ either $B \lq A$ or
$B \lq A^*$.   Let $u$ be the set of those $A \in \E$ for which $B \lq A^*$ but $A^*$
is minimal with this property.  Then if $A', A  \in u$ and $A' \subset A$, and $A' \lq C \lq A^*$
for $C \in \E$, then if $B \lq C^*$ then $A' = C$ while if $B \lq C$ then $A = C$ by the minimality
of $A^*$.  Thus $u \in VT$.  
\end{proof}
\medskip
The above decomposition of a graph is well-known.  Tutte \cite{[Tu]} extended this approach, 
by analyzing how finite graphs can be separated by removing pairs of points. He 
obtained a tree decomposition in which some vertices corresponded to
  $3$-blocks or maximal $3$-connected subgraphs and Droms, Servatius and Servatius \cite{[DrSS]} generalized Tutte's result to infinite locally finite graphs.
We give a new proof of this latter result.

Let $X$ be a locally finite connected graph which is $2$-connected but not $3$-connected.
 This means that there are no cut-points, i.e. vertices whose removal disconnects the graph,
 but $X$ can be disconnected by the removal of a pair of vertices.
 We also assume that $X$ is simple (there is at most one edge joining any pair of vertices and there are no loops) and
there are no vertices of degree $2$.  If a graph does not have this property then, unless the graph
is $2$-regular, or contains a subdivided loop,  it can be replaced by a homeomorphic graph in which there are no vertices 
of degree $2$ and if there are multiple edges joining two vertices  they can be replaced by
a single edge joining those vertices.

 We consider subgraphs  $A \subset X$ with the following property
 \begin{itemize}
 \item [{(i)}] there are exactly $2$ vertices $u,v \in VA$ such that if $w \in VA, w \not= u, w\not=v$
 then every edge of $X$ incident with $w$ is also in $A$. At least one, but not all, of the edges incident with
 each of $u, v$ are in $A$. 
 
 \item[{(ii)}] $A$ does not consist of a single edge  and it does not consist of a subgraph obtained
 from $X$ by removing one edge and no vertices.
\end{itemize}
The set of such subgraphs is denoted ${\cb}_2X$.  For $A$ as above we put $\delta A = \{ u,v \}$.
Following Tutte \cite{[Tu]}, we call this the {\it hinge} of $A$.
Since $X$ is connected every vertex in $A$ is joined to either $u$ or $v$ by a path.
Thus $A$ has at most $2$ components.  But if there were $2$ components then
$u, v$ would be cut-points which would contradict our assumption that 
$X$ is $2$-connected.  It follows that $A$ is connected.
 \begin{lemm}\label{lemm:sep} Let $u \in VX$.  There are only finitely many $A \in {\cb}_2X$ such that
$u \in \delta A$.
\end{lemm}
\begin{proof} Let $\delta A = \{ u, w \}$ and let $v \in VX$ be adjacent to $u$.
If $v \in VA$ then there is another vertex $p \in VA^*$ which is adjacent to $u$.  This is because 
some but not all edges incident with $u$ are in $A$.  Similarly if $v \in VA^*$, then there is a vertex
$p$ adjacent to $u$ which is in $VA$.  For each $p$ adjacent to $u, p \ne v$ choose a path
in $X$ joining $v$ and $p$ not passing through $u$.  Such a path exists since $u$ is not a cut-point.
There are finitely many of these paths since $u$ has finite degree.  As one of these paths joins 
vertices in $A$ and $A^*$ it must pass through $w$.  We see then that there are only
finitely many possibilities for $w$.  Now $A$ is determined by the pair $\delta A = \{ u, w\}$ and the edges
of $A$ which are incident with this pair.   There are only finitely many possibilities and so the
lemma is proved.
\end{proof}

\begin{lemm}\label{lemm:2.4}  Let $x \in VX$ and let $A_1 \supseteq A_2 \supseteq \dots $
be a sequence of subgraphs in ${\cb}_2X$ each containing $x$. 
The sequence is eventually constant.
\end{lemm}
\noindent
\begin{proof} Let $v \in VA_1^*$ and let $p$ be a path in $X$ joining $v, x$.  Now 
for each $i$, the path $p$ joins vertices in $A_i $ and $A_i^*$ and so it must intersect $\delta A$.  As the path $p$ contains a fixed finite set of vertices the Lemma
follows from Lemma ~\ref{lemm:sep} .
\end{proof}
\begin{lemm}\label{lemm:2.5}  Let $x \in VX$ and let $A $ be a smallest element of
${\cb}_2X$ containing $x$.   Let $B \in {\cb}_2X$ then at least one of
$A \subseteq B,  A \subseteq  B^*, A^* \subseteq B, A^* \subseteq B^*$.
\end{lemm}
\noindent
\begin{proof} Let  $\delta A = \{ u, v\} , \delta B = \{ p, q \}$.
Suppose first that $u, v, p, q$ are all distinct.   If $A$ contained both $p$ and $q$
and $u, v \in B$ then every vertex in $A^* \cap B^*$ would have every incident
edge in  $A^* \cap B^*$ and so it is $\emptyset $ since it cannot be $X$. Hence $A^* \subseteq B$.
Let $p,q, u \in A \cap B, v \in B^*$ then every vertex in  $A^* \cap B^*$ except $v$ would 
have every incident edge in  $A^* \cap B^*$.  Since $v$ is not a cut-point,  $A^* \cap B^*= \{ v \}$.
Similarly $A^*\cap B = \{ u \}$ which is absurd, since $A^*$  does not consist of just two vertices.
So this case does not occur.
The case when the hinge  for $A$ is contained in $B$ or $B^*$ and only one vertex  of  $\delta B$ is contained in $A$ cannot occur either,
by reversing the roles of $A$ and $B$ in the above argument.
Still assuming that $u,v,p,q$ are all distinct, suppose $x, u,  p \in A\cap B$,
and  $v \in  B^*, q \in A^* $.
Now $A \cap B$ would be a smaller
element of ${\cb}_2X$ than $A$ containing $x$, which is a contradiction unless $A\cap B$ consists of a single edge $e$ with   two vertices $u, p$ and $x = u$ or $x = p$.  If $x=u$, then removing  the edge 
$e$ and vertex $p$ from $A$ produces a smaller element of ${\cb}_2X$ unless $A$ has $2$ edges and $3$ vertices
$u, p, v$. Similarly if $x = p$ then removing $e$ and the vertex $u$ from $A$ produces a  smaller element of ${\cb}_2X$ unless $A$ has $2$ edges and $3$ vertices   $u, p, v$.
In both these cases $p$ has degree  $2$ which contradicts our hypothesis.

 All 
other cases for $\{ u,v, p,q \} $ all distinct can be reduced to one of the cases considered
by a relabelling.  Thus if  $u,v,p,q$ are all distinct, then one of
the four inclusions occurs.

Suppose now that $u = p$ and $u, v, q $ are distinct. If $A$ contains $q$ and $B$ contains 
$v$ then every vertex of $A^*\cap B^*$ except $u$ would have every incident edge
in  $A^*\cap B^*$ which means that $A^*\cap B^* = \{ u \}$.  But then $A^* \subseteq B\cup \{ u \} = B$.

Finally if $u = p, v = q$ and $x \in B$ then either $A\cap B$ contains  a single edge $e$ joining vertices $u , v$ one of which is $x$  or $A\cap B$ is an element of 
${\cb}_2X $ containing $x$ and so  by the minimality of $A, A= A\cap B$ and $A\subseteq B$.
 In the case when $A\cap B$ contains  a single edge $e$ joining vertices $u , v$ one of which is $x$, 
 note that $x \in B^*$ and $A\cap B^*$ must contain more than $2$ vertices, and so by the minimality
 of $A, A = A\cap B^*$.  
\end{proof}
\medskip

Choose $A_0$ to be a fixed smallest element of ${\cb}_2X$ containing a particular $x_0$.
If we choose $\E$ to be the set $\{ gA_0, gA_0^* | g \in G \} $ then we have
the following theorem.
\begin{theo}\label{theo:2.6}  Let $X$ be a $G$-graph which is locally finite and $2$-connected but
not $3$-connected.  Suppose also that $X$ is not a cycle, that it is simple graph and has at least $4$ vertices. Then there is a non-empty $G$-subset $\E \subset {\cb}_2$
closed under the order reversing
involution $A \mapsto A^*$, satisfying the following two conditions
\item {(i)} for every $A, B \in \E$ at least  one of $A \lq B, A\lq B^*, B \lq A, B \lq A^*$ holds.
\item {(ii)} if $A, B \in \E$ and $A \lq B$ then there are only finitely many $C \in E$ such
that $A \lq C \lq B$.  
\end{theo}

\begin{proof} It remains to prove  that $\E$ satisfies (ii).  Let $A, B \in \E, A\lq B$.  Let $u \in VA, w \in B^*$ and 
let $u =v_1, v_2, \dots , v_n = w$ be the vertices of a path joining them.
If $A \lq C \lq B$ then $v_i \in \delta C$ for some $i = 1, 2, \dots , n$.
It follows easily from  Lemma~\ref{lemm:sep}  that $\E$ satisfies (ii).
\end{proof}

As above the $G$-set $(\E , \lq )$ satisfying the above conditions can be regarded as the oriented edge set of a $G$-tree $T = T(\E )$. 

Also as above a vertex $v \in VT$ is an equivalence class under the equivalence relation $\sim $ on $\E $ defined as follows:-

$A\sim B$ if  $A \subset  B^*$ but for no $C \in \E$ is $A\subset C \subset B^*$.

As in the analysis of ${\cb}_1X$ some vertices correspond to hinges.
There will be other vertices $v \in VT$ unless for every $A \in \E$ there is only one other
$B \in \E$ (namely $A^*$) such that $\delta B = \delta A$.

For any such $v$ we now show how to associate a graph $Z = Z_v$.
Let $Z' = \bigcap _{E\in v}E^*$.  
Clearly if $A \lq  B^*$ then $ B^*$ contains the hinges of both $A$ and $B$.
 It follows easily that $Z'$ will contain every vertex which is in $\delta E $
for any $E \in v$. In particular it is non-empty.  
The graph $Z$ is obtained from $Z'$ by adjoining an extra edge for each $E \in v$
 joining the  vertices of the hinge of $E$ provided these vertices are not already
 joined by an edge..  In [Tu] these edges are called virtual edges.
 Note that for each such extra edge 
 (corresponding to $E$) there is a path in $E$ joining the vertices of the edge.
 If we let $Z''$ be the graph consisting of $Z'$ together with each such path 
 then $Z''$ is isomorphic to a subdivision of $Z$ and it is a subgraph of $X$. Thus
 if $X$ is planar, then so is each $Z_v$.
 Let $G_v$ be the stabilizer of $v$, then $Z$ is a $G_v$-graph.
 
 We now show that if for any $v \in VT$, not corresponding to a hinge, $Z_v$ is not $3$-connected and it is not a cycle, then we can 
 enlarge the $G$-set $\E$.
 
\begin{step}  The graph $Z$ is connected.
 \end{step}
 Let $x = x_1, x_2, \dots , x_n = y$ be the vertices of a shortest path in $X$ joining two
 vertices of $VZ = VZ'$.  Any part of the path which lies in a particular $E \in v$ must start and
 end with a vertex which is special for $E$.  Since we have a shortest path, these vertices are different and so joined by an edge in $Z$. We obtain a path in $Z$ by replacing each
 such section by the corresponding edge in $Z$.   It follows that $Z$ is connected.
 
 \begin{step}  The graph $Z$ is $2$-connected.
 \end{step}
 Suppose $x$ is a vertex whose removal disconnects $Z$.  Let $U$ be one of the components
 of $Z - \{ x \}$.  We want to show that $x$ is a cut-point for $X$ as well as for $Z$.
 Let $U'$ be the subgraph of $X$ consisting of those vertices and edges which are joined
 to a vertex of $U$ by a path in $X - \{ x \}$.  Suppose $z \in VZ \cap U'$.  There is a path in
 $X - {x}$ joining $z$ to a vertex of $U$.  Now as in Step 1 replace this path by a path
 in $Z$.  This path will still not pass through $x$ as the vertices are a subset of the original path.
 But $U$ is a component of $Z - \{ x \}$ and so the path must lie in $U$.  Hence $VZ \cap U' = U$
 and $x$ is a cut-point for $X$.
 
  \begin{step} If $Z$ is not $3$-connected and it is not a cycle (i.e. it contains
 vertices of degree greater than $2$), then there is non-empty $G_v$-set
 $\E _v \lq {\cb}_2Z$ satisfying conditions (i)  above.
 \end{step}
 This follows immediately from Theorem~\ref{theo:2.6}.
 
 \begin{step}  If $Z$ is not $3$-connected and it is not a cycle we can enlarge $\E$, i.e.
 there is a $G$-set $\E ' \lq {\cb}_2X$ properly containing $\E$ and satisfying conditions
 (i) and (ii).
 \end{step}
 Let $A_v \in \E _v$ with  $\delta A_v = \{ x,y \} $ and let $A_v'$ be the subgraph of $X$ 
 consisting of edges and vertices which are joined to $A_v$ be a path lying in $X - \{ x, y\} $.
 As in the proof of Step 2, $A_v' \cap Z_v = A_v$. Also if $A_v, B_v \in \E _v$ and $A_v \lq B_v$
 then $A_v' \lq B_v'$.
 Let $\E' = \E \cup \{ gA_v', gA_v'^* | g \in G, A_v \in \E _v\} $.  In order to show that $\E '$ satisfies (iii)
 it suffices to show that if $g \in G$ then $A_v', gB_v'$ are nested.  We know this is the case
 if $g \in G_v$.  If $g\notin G_v$ then $v \not= gv$ and $gB_v \in \E _{gv}$.
 Let the path in $T$ joining $v$ and $gv$ start with edge $e$ and finish with edge $f$, so that
 $\iota e = v, \tau f = gv$.  Thus $e$ is in the equivalence class defining $v$ and $f^*$ is in the equivalence class defining $gv$.  Also because $e, f$ are oriented coherently as part of
 a edge path, we have $f \lq e$.  Now $e$ corresponds to the subgraph of $Z _v$ consisting of 
 a single virtual edge and its vertices. Hence either $e \lq A_v$ or $e \lq A_v^*$.  Similarly  $f^* \lq gB_v'$ or
 $f^* \lq gB_v'*$.  But if, say, $e\lq A_v'$ and $f^*\lq gB_v'$, then $A_v'^* \cap gB_v'^*
 \lq e^*\cap f = \emptyset $ and so $A_v'$ and $gB_v'$ are nested.
 
 \medskip
 It follows from the previous discussion  that if we can show that we can bound the number of times 
 the $G$-set $\E $ can be enlarged, then there must be a set $E$ for which all the corresponding
 vertex graphs $Z_v$ are $3$-connected. 
We show that this is the case if $X$ is $G$-finite. In Lemma~\ref{lemm:sep}, it was shown that
there are only finitely many $A \in {\cb}_2X$ for which $\delta A$ contains a particular vertex.  It follows easily that the size of the $G$-set $\E$ is bounded and by choosing a largest such $\E$ the corresponding $Z_v$'s  are $3$-connected. These are the $3$-blocks of $X$.
Note that we have to allow $Z_v$ to have $2$ or $3$ vertices, which we regard as $3$-connected.

We can now combine our results so that we consider all graphs which are not 
$3$-connected.  We still assume that $X$ is a simple graph with no vertices of degree $2$.
Enlarge ${\cb}_2X$ to consist of subgraphs $A$ for which either
\begin{itemize} 
\item [(i)] either there are $2$ vertices $u,v \in VA$ such that if $w \in VA, w \not= u, w\not=v$
 then every edge of $X$ incident with $w$ is also in $A$.
 
 or 
 
 \item[(ii)] there  is   a single vertex $u \in  VA$ such that if $w \in VA, w \not= u$,
 then every edge of $X$ incident with $w$ is also in $A$.
 
 If $X$ is $2$-connected, there are no subgraphs satisfying (ii).
\end{itemize}

\begin{theo}\label{theo:2.7}  Let $X$ be a $G$-finite locally finite graph, where $G$ is the
automorphism group of $X$.   There is a $G$-finite $G$-subset $\E$ of ${\cb}_2X$
which satisfies the conditions (i) and (ii).  If $T$ is the $G$-tree associated with
$\E$ then every edge of $T$ is associated with a decomposition $(A, A^*)$
where $A, A^*$ are subgraphs such that $A\cup A^* = X$ and $A\cap A^*$ consists of
either one or two vertices.  To every vertex $v \in VT$ corresponds either a cut-point of $X$, a $G_v$- cycle or a $3$-connected $G_v$-graph
which has  a subdivision which is a subgraph of  $X$, so that it is planar if $X$ is planar.

 \end{theo}

\bigskip

In \cite{[Tu]} Tutte gives examples where the graph associated with a vertex is a cycle.

\bigskip

Theorem~\ref{theo:2.7}  and Theorem~\ref{theo:main} give us a lot of information about $G$-finite locally finite planar graphs.
In particular if the $G$-action is free, as it is for a Cayley graph, then we can remove
the $3$-connected condition from Theorem~\ref{theo:main}.
\begin{theo}\label{theo:2.8} If $G$ is a group and $X$ is  a  connected locally finite planar graph $X$ on
which $X$ acts freely  so that
$G\backslash X$ finite, then $G$ (or an index two subgroup of $G$) is the fundamental group
of a graph of groups in which each vertex group is either a planar discontinuous group
or a free product of finitely many  cyclic groups and all edge groups are finite cyclic groups (possibly
trivial).  
\end{theo}
\noindent 
\begin{proof} If the action is free, then the stabilizer of any cut-point is trivial and the stabilizer
of any edge of the tree $T$  in Theorem~\ref{theo:2.7} has order at most $2$.  This means that $G$
is the fundamental group of a graph of groups in which each vertex group $G_v$  acts freely on
a planar $3$-connected planar graph $X_v$ and $G_v\backslash X_v$ is finite and
each edge group has order at most $2$.  The result follows immediately from
Theorem~\ref{theo:2.7}.
\end{proof}

\end{document}